# Asymptotics of solving a singularly perturbed system of transport equations with fast and slow components in the critical case[*]

**Nesterov A.V.** [1][0000-0002-4702-4777]

[1] PLEKHANOV Russian University of Economics, Stremyanny lane 36, Moscow, 117997, Russia
`andrenesterov@yandex.ru`
[2]

**Abstract.** An asymptotic expansion with respect to a small parameter of the solution of the Cauchy problem is constructed for a system of three transport equations, two of which are singularly perturbed with degeneracy of the entire high-order part of the transfer operator, and the third equation does not explicitly contain a small parameter. A special feature of the problem is that it belongs to the so-called critical case: the solution of the degenerate problem is a one-parameter family. The asymptotic expansion of the solution under smooth initial conditions is constructed as a sum of the regular part and boundary functions. It is noted that in the case of discontinuous initial conditions, the regular part will have a gap on the "pseudo-characteristics" of the initial system that do not coincide with the characteristics themselves.

## 1 Introduction

Singularly perturbed systems of hyperbolic equations and transport equations have been considered in many papers, and some of them are discussed in [1]. In [2], the author studied the asymptotics of the small parameter solution of an initial-boundary value problem for a hyperbolic system of two singularly perturbed equations in the critical case [3] ( the solution of the degenerate system is not unique)

$$\begin{cases} \mu(u_t + \Lambda_1 u_x) = -au + bv, \\ \mu(v_t + \Lambda_2 v_x) = \gamma(au - bv), x > 0, t > 0 \end{cases}$$

with initial-boundary conditions

---

[*] This research was performed in the framework of the state task in the field of scientific activity of the Ministry of Science and Higher Education of the Russian Federation, project "Models, methods, and algorithms of artificial intelligence in the problems of economics for the analysis and style transfer of multidimensional datasets, time series forecasting, and recommendation systems design", grant no. FSSW-2023-0004.

$$\begin{cases} u(x,0) = \varphi_1(x), v(x,0) = \varphi_2(x), \\ u(0,t) = \psi_1(t), v(0,t) = \psi_2(t), \end{cases}$$

where $\Lambda_2 > \Lambda_1 > 0, a > 0, \gamma b > 0, 0 < \mu \ll 1$. When $\mu = 0$, solution is a one-parameter family $u(x,t,0) = \phi(x,t), v(x,t,0) = \dfrac{a}{b}\phi(x,t)$, where $\phi(x,t)$ is an arbitrary smooth function. In [2], it is shown that $\phi(x,t)$ there is a solution of the initial-boundary value problem for the transport equation $\phi_t + \Lambda_3 \phi_x = 0, x > 0, t > 0, \Lambda_2 > \Lambda_3 > \Lambda_1 > 0$ with conditions inconsistent at the corner point of the boundary $(0,0)$ and has a discontinuity of the first kind on a line $x = \Lambda_3 t$, that is not a characteristic of the original system of equations. In [4], it is shown that the exact solution has the character of an asymptotically narrow transition layer in the vicinity of this line, and an asymptotic expansion of the solution in the vicinity of the transition layer zone is obtained. In [5]-[6], asymptotics of solutions to more general singularly perturbed systems of transport equations were constructed, but in these works, systems of equations containing only fast components were studied, i.e. all equations contained small parameters with the highest part. It is of interest to study similar systems in which part of the equations contains small parameters with higher terms of operators (fast components), and part - does not contain (slow components). Singularly perturbed systems of ordinary differential equations containing both fast and slow components are studied in detail in [3].

## 2    Problem statement

Consider a singularly perturbed system of transport equations

$$\begin{cases} \varepsilon^2(u_t + k_1 u_x) = -au + bv + \varepsilon^2 c_1 w, \\ \varepsilon^2(v_t + k_2 v_x) = au - bv + \varepsilon^2 c_2 w, \\ w_t + k_3 w_x = a_3 u + b_3 v + c_3 w, \end{cases} \qquad (1)$$

with initial conditions

$$\begin{cases} u(x,0,\varepsilon) = u^0(x), \\ v(x,0,\varepsilon) = v^0(x), \\ w(x,0,\varepsilon) = w^0(x). \end{cases} \qquad (2)$$

where: $\{u(x,t,\varepsilon), v(x,t,\varepsilon), w(x,t,\varepsilon)\}$ - solution (below the entry, the explicit dependence on the parameter $\varepsilon$ can be omitted), $\{x,t\} \in H = \{|x| < \infty; 0 \le t \le T, T > 0\}; 0 < \varepsilon \ll 1$ - small parameter, , the functions $u^0(x), v^0(x), w^0(x)$ have a certain smoothness (as discussed below) and are finite or rapidly decreasing $|u^0(x)|, |v^0(x)|, |w^0(x)| < Ce^{-\kappa x^2}, C, \kappa > 0$.

The singularly perturbed system (1) contains both fast ($u(x,t,\varepsilon), v(x,t,\varepsilon)$) and slow $w(x,t,\varepsilon)$ components. In addition, it belongs to the so-called critical case [3], because the solution of the degenerate problem is a one-parameter family. If $\varepsilon = 0$, the task takes the form $-au + bv = 0, w_t + k_3 w_x = a_3 u + b_3 v + c_3 w$. Hence $v(x,t,0) = \frac{a}{b} u(x,t,0)$, where $u(x,t,0)$ is an arbitrary smooth function. The function $w(x,t,0)$ is defined from the equation $w_t + k_3 w_x = (a_3 + \frac{ab_3}{b})u + c_3 w$. Systems of this type containing only fast components have been studied in a number of papers [2] - [6] and others. The presence of slow components significantly changes the structure of the asymptotic behavior of the solution and, consequently, the behavior of the solution itself.

## 3  Construction of an AE function under smooth initial conditions.

And the asymptotic expansion (AE) with respect to a small parameter of the solution of problem (1)-(2), if the smoothness condition is satisfied on the initial data, under the conditions of smoothness of the initial conditions, we will look for in the form

$$\begin{pmatrix} u(x,t,\varepsilon) \\ v(x,t,\varepsilon) \\ w(x,t,\varepsilon) \end{pmatrix} = \overline{Y} + \Pi = \begin{pmatrix} \overline{u}(x,t,\varepsilon) \\ \overline{v}(x,t,\varepsilon) \\ \overline{w}(x,t,\varepsilon) \end{pmatrix} + \begin{pmatrix} \Pi u(x,\tau,\varepsilon) \\ \Pi v(x,\tau,\varepsilon) \\ \Pi w(x,\tau,\varepsilon) \end{pmatrix} = $$
$$= \sum_{i=0}^{N} \varepsilon^i \left( \begin{pmatrix} \overline{u}_i(x,t) \\ \overline{v}_i(x,t) \\ \overline{w}_i(x,t) \end{pmatrix} + \begin{pmatrix} \pi_i u(x,\tau) \\ \pi_i v(x,\tau) \\ \pi_i w(x,\tau) \end{pmatrix} \right) + \begin{pmatrix} R_N u(x,t,\varepsilon) \\ R_N v(x,t,\varepsilon) \\ R_N w(x,t,\varepsilon) \end{pmatrix} = \qquad (3)$$
$$= \overline{Y}_N + \Pi_N + R_N$$

where $\overline{u}, \overline{v}, \overline{w}$ is the regular part of AP, $\Pi u, \Pi v, \Pi w$ - boundary functions, $R$ is the residual term, $\tau = \frac{t}{\varepsilon^2}$. In this paper, we restrict ourselves to a detailed description of the construction only the main members of AE (3).

In the case of smooth initial conditions, the construction of AE (3) does not cause any special difficulties and is carried out by the method of boundary functions of A. B. Vasilyeva and V. F. Butuzov [3]. Of greater interest is the construction of the AE in the case when the initial conditions are piecewise continuous or have regions of large gradients.

### 3.1. Construction a regular part of the AE

For smooth initial conditions, the construction of a regular part of the AE is not difficult, but it is important for a more interesting case-initial conditions of the "burst" type or with discontinuities of the first kind.

The regular part of AE (3) has the form

$$\begin{pmatrix} \bar{u}(x,t,\varepsilon) \\ \bar{v}(x,t,\varepsilon) \\ \bar{w}(x,t,\varepsilon) \end{pmatrix} = \sum_{i=0}^{\infty} \varepsilon^i \begin{pmatrix} \bar{u}_i(x,t) \\ \bar{v}_i(x,t) \\ \bar{w}_i(x,t) \end{pmatrix} \qquad (4)$$

Following [3], to determine the terms of expansion (4), we substitute expansion (4) in the system (1)

$$\begin{cases} \varepsilon^2(\bar{u}_{0,t} + k_1\bar{u}_{0,x}) + \ldots = -a\bar{u}_0 + b\bar{v}_0 + \varepsilon(-a\bar{u}_1 + b_1\bar{v}_1) + \\ \quad + \varepsilon^2(-a\bar{u}_2 + b\bar{v}_2 + c_1\bar{w}_0) + \ldots, \\ \varepsilon^2(\bar{v}_{0,t} + k_2\bar{v}_{0,x}) + \ldots = a\bar{u}_0 - b\bar{v}_0 + \varepsilon(a\bar{u}_1 - b_1\bar{v}_1) + \\ \quad + \varepsilon^2(a\bar{u}_2 - b\bar{v}_2 + c_2\bar{w}_0) + \ldots, \\ \bar{w}_{0,t} + k_3\bar{w}_{0,x} + \varepsilon(\bar{w}_{1,t} + k_3\bar{w}_{1,x}) + \ldots = a_3\bar{u}_0 + b_3\bar{v}_0 + c_3\bar{w}_0 + \\ \quad + \varepsilon(a_3\bar{u}_1 + b_3\bar{v}_1 + c_3\bar{w}_1) + \ldots. \end{cases}$$

and we equate the left and right terms with the same powers of the parameter $\varepsilon$.

If the parameter degree is zero, we get

$$\varepsilon^0 : \begin{cases} -a\bar{u}_0 + b\bar{v}_0 = 0, \\ a\bar{u}_0 - b\bar{v}_0 = 0, \\ \bar{w}_{0,t} + k_3\bar{w}_{0,x} = a_3\bar{u}_0 + b_3\bar{v}_0 + c_3\bar{w}_0. \end{cases} \qquad (5)$$

From the first two equations

$$\bar{v}_0 = \frac{a}{b}\bar{u}_0 \qquad (6)$$

where $\bar{u}_0(x,t)$ - is an undefined function.

For the first degree of the parameter, we get

$$\varepsilon^1 : \begin{cases} -a\bar{u}_1 + b\bar{v}_1 = 0, \\ a\bar{u}_1 - b\bar{v}_1 = 0, \\ \bar{w}_{1,t} + k_3\bar{w}_{1,x} = a_3\bar{u}_1 + b_3\bar{v}_1 + c_3\bar{w}_1 \end{cases}$$

From the first two equations

$$\bar{v}_1 = \frac{a}{b}\bar{u}_1, \quad \text{where } \bar{u}_1(x,t) \text{ - is an undefined function.}$$

For the second power of the parameter, we get

$$\begin{cases} \bar{u}_{0,t} + k_1\bar{u}_{0,x} = -a\bar{u}_2 + b\bar{v}_2 + c_1\bar{w}_0, \\ \bar{v}_{0,t} + k_2\bar{v}_{0,x} = a\bar{u}_2 - b\bar{v}_2 + c_2\bar{w}_0, \\ \bar{w}_{2,t} + k_3\bar{w}_{2,x} = a_3\bar{u}_2 + b_3\bar{v}_2 + c_3\bar{w}_2 \end{cases}$$

The first two equations are a system of two linear equations with respect to $\bar{u}_2, \bar{v}_2$. The solvability condition of this system gives

$$\bar{u}_{0,t} + \bar{v}_{0,t} + k_1\bar{u}_{0,x} + k_2\bar{v}_{0,x} = (c_1+c_2)\bar{w}_0 \qquad (7)$$

Substituting relation (6) in (7), we obtain the equation for determining the function $\bar{u}_0$

$$\left(1+\frac{a}{b}\right)\bar{u}_{0,t} + \left(k_1+\frac{a}{b}k_2\right)\bar{u}_{0,x} = (c_1+c_2)\bar{w}_0,$$

then, in conjunction with the third equation (5), gives a system of equations for determining $\bar{u}_0, \bar{w}_0$. Let's introduce the notation

$$\frac{bk_1+ak_2}{a+b}=V, \frac{b(c_1+c_2)}{a+b}=C, \frac{a_3b+ab_3}{b}=c, D=Cc=\frac{(c_1+c_2)(a_3b+ab_3)}{a+b}$$

Using these notations, the system of equations for the zero terms of the expansion takes the form

$$\begin{cases} \bar{u}_{0,t} + V\bar{u}_{0,x} = C\bar{w}_0, \\ \bar{w}_{0,t} + k_3\bar{w}_{0,x} = c\bar{u}_0 + c_3\bar{w}_0. \end{cases} \qquad (8)$$

System (8) can easily be reduced to a single second-order hyperbolic equation. Differentiating the first equation with respect to t, we obtain

$$\bar{u}_{0,tt} + V\bar{u}_{0,xt} = C\bar{w}_{0,t},$$

Differentiating the first equation with respect to x, we obtain

$$\bar{u}_{0,tx} + V\bar{u}_{0,xx} = C\bar{w}_{0,x},$$

Substituting these relations in (8), we obtain

$$\bar{w}_{0,t} + k_3\bar{w}_{0,x} = \frac{1}{C}(\bar{u}_{0,tt} + V\bar{u}_{0,xt} + k_3(\bar{u}_{0,tx} + V\bar{u}_{0,xx})) = c\bar{u}_0 + c_3\bar{w}_0 =$$

$$= c\bar{u}_0 + c_3\frac{1}{C}(\bar{u}_{0,t} + V\bar{u}_{0,x})$$

Excluding $\bar{w}_0$, we obtain the equation for $\bar{u}_0$

$$\bar{u}_{0,tt} + (V + k_3)\bar{u}_{0,xt} + k_3 V \bar{u}_{0,xx} = D\bar{u}_0 + c_3(\bar{u}_{0,t} + V\bar{u}_{0,x}) \tag{9}$$

This equation can also be written as

$$(\frac{\partial}{\partial t} + V\frac{\partial}{\partial x})(\frac{\partial}{\partial t} + k_3\frac{\partial}{\partial x})\bar{u}_0 = c_3(\frac{\partial}{\partial t} + V\frac{\partial}{\partial x})\bar{u}_0 + D\bar{u}_0 \tag{10}$$

For the remaining terms of the expansion, similar inhomogeneous equations are obtained. Other expansion terms (4) they are constructed in exactly the same way [3].

The coefficient $V$ satisfies the inequality $\min(k_1, k_2) < V < \max(k_1, k_2)$, so the characteristics of the equation (10) - line $x - Vt = const, x - k_3 t = const$ generally do not coincide with the characteristics of the original system (1) $x - k_1 t = const, x - k_2 t = const, x - k_3 t = const$.

*Remark.* In accordance with the construction algorithm

$$\bar{Y}_N = \sum_{i=0}^{N} \varepsilon^i \begin{pmatrix} \bar{u}_i(x,t) \\ \bar{v}_i(x,t) \\ \bar{w}_i(x,t) \end{pmatrix} \equiv \begin{pmatrix} \bar{U}_N(x,t) \\ \bar{V}_N(x,t) \\ \bar{W}_N(x,t) \end{pmatrix} \tag{11}$$

where N is a natural number, satisfies the system

$$\begin{cases} \varepsilon^2(\bar{U}_{N,t} + k_1\bar{U}_{N,x}) = -a\bar{U}_N + b\bar{V}_N + \varepsilon^2 c_1\bar{W}_N + O(\varepsilon^{N+1}), \\ \varepsilon^2(\bar{V}_{N,t} + k_2\bar{V}_{N,x}) = a\bar{U}_N u - b\bar{V}_N + \varepsilon^2 c_2\bar{W}_N + O(\varepsilon^{N+1}), \\ \bar{W}_{N,t} + k_3\bar{W}_{N,x} = a_3\bar{U}_N + b_3\bar{V}_N + c_3\bar{W}_N + O(\varepsilon^{N+1}). \end{cases} \tag{12}$$

### 3.2 Construction of boundary functions.

The solution of the system (8) containing two equations cannot satisfy the three initial conditions (2). In order to satisfy the initial conditions, we construct boundary functions following [3]

$$\begin{pmatrix} \Pi u(x,\tau,\varepsilon) \\ \Pi v(x,\tau,\varepsilon) \\ \Pi w(x,\tau,\varepsilon) \end{pmatrix} = \sum_{i=0}^{\infty} \varepsilon^i \begin{pmatrix} \pi_i u(x,\tau) \\ \pi_i v(x,\tau) \\ \pi_i w(x,\tau) \end{pmatrix}, \tau = t/\varepsilon^2 \tag{13}$$

Boundary functions formally satisfy the system

$$\begin{cases} \varepsilon^2(\Pi u_t + k_1 \Pi u_x) = -a\Pi u + b\Pi v + \varepsilon^2 c_1 \Pi w, \\ \varepsilon^2(\Pi v_t + k_2 \Pi v_x) = a\Pi u - b\Pi v + \varepsilon^2 c_2 \Pi w, \\ \Pi w_t + k_3 \Pi w_x = a_3 \Pi u + b_3 \Pi v + c_3 \Pi w, \end{cases} \tag{14}$$

together with the regular part, they meet the initial conditions

$$\begin{cases} \bar{u}(x,0)+\Pi u(x,0)=\overset{0}{u}(x) \\ \bar{v}(x,0)+\Pi v(x,0)=\overset{0}{v}(x), \\ \bar{w}(x,0)+\Pi w(x,0)=\overset{0}{w}(x). \end{cases} \qquad (15)$$

and are functions of the boundary layer:

$$\begin{cases} \Pi u(x,\tau) \to 0, \tau \to +\infty, \\ \Pi v(x,\tau) \to 0, \tau \to +\infty, \\ \Pi w(x,\tau) \to 0, \tau \to +\infty. \end{cases} \qquad (16)$$

Passing to variables $x, \tau = t/\varepsilon^2$ in system (14), we obtain

$$\begin{cases} \Pi u_\tau + \varepsilon^2 k_1 \Pi u_x = -a\Pi u + b\Pi v + \varepsilon^2 c_1 \Pi w, \\ \Pi v_\tau + \varepsilon^2 k_2 \Pi v_x = a\Pi u - b\Pi v + \varepsilon^2 c_2 \Pi w, \\ \Pi w_\tau + \varepsilon^2 k_3 \Pi w_x = \varepsilon^2 a_3 \Pi u + \varepsilon^2 b_3 \Pi v + \varepsilon^2 c_3 \Pi w, \end{cases} \qquad (17)$$

Substituting the expansion (13) into the system (17), we obtain equations for the terms of the expansion (13). In the zero approximation

$$\begin{cases} \pi_0 u_\tau = -a\pi_0 u + b\pi_0 v, \\ \pi_0 v_\tau = a\pi_0 u - b\pi_0 v, \\ \pi_0 w_\tau = 0. \end{cases} \qquad (18)$$

From equations (18), taking into account the conditions

$$\begin{cases} \pi_0 u(x,\tau) \to 0, \tau \to +\infty, \\ \pi_0 v(x,\tau) \to 0, \tau \to +\infty, \\ \pi_0 w(x,\tau) \to 0, \tau \to +\infty \end{cases}$$

should

$$\begin{cases} \pi_0 u(x,\tau) = \psi(x)e^{-(a+b)\tau}, \\ \pi_0 v(x,\tau) = -\psi(x)e^{-(a+b)\tau}, \\ \pi_0 w = 0, \end{cases} \qquad (19)$$

where $\psi(x)$ is an arbitrary function.

Substituting (6) and (19) into the initial conditions (15), we obtain

$$\begin{cases} \bar{u}(x,0) + \psi(x) = u^0(x) \\ \dfrac{a}{b}\bar{u}(x,0) - \psi(x) = v^0(x), \\ \bar{w}(x,0) = w^0(x). \end{cases}$$

Hence we obtain the initial conditions for the functions $\bar{u}(x,0), \bar{w}(x,0)$

$$\begin{cases} \bar{u}(x,0) = \dfrac{b}{a+b}(u^0(x) + v^0(x)), \\ \bar{w}(x,0) = w^0(x) \end{cases} \quad (20)$$

as well as the function $\psi(x)$

$$\psi(x) = \dfrac{au^0(x) - bv^0(x)}{a+b} \quad (21)$$

The remaining terms of the expansion (13) are constructed in exactly the same way as [3] and satisfy the bound

$$|\pi_i u(x,\tau)|, |\pi_i v(x,\tau)|, |\pi_i w(x,\tau)| < C\exp(-\kappa\tau), C, \kappa > 0 \quad (22)$$

Constants $C, \kappa$ generally depend on the number $i$. The proof of estimates (22) repeats almost verbatim the proof in [3].

*Remark.* In accordance with the construction algorithm

$$\begin{pmatrix} \Pi_N u(x,\tau) \\ \Pi_N v(x,\tau) \\ \Pi_N w(x,\tau) \end{pmatrix} = \sum_{i=0}^{N} \varepsilon^i \begin{pmatrix} \pi_i u(x,\tau) \\ \pi_i v(x,\tau) \\ \pi_i w(x,\tau) \end{pmatrix} \quad (23)$$

satisfies the system

$$\begin{cases} \Pi_N u_\tau + \varepsilon^2 k_1 \Pi_N u_x = -a\Pi_N u + b\Pi_N v + \varepsilon^2 c_1 \Pi_N w + O(\varepsilon^{N+1} e^{-\kappa\tau}), \\ \Pi_N v_\tau + \varepsilon^2 k_2 \Pi_N v_x = a\Pi_N u - b\Pi_N v + \varepsilon^2 c_2 \Pi_N w + O(\varepsilon^{N+1} e^{-\kappa\tau}), \\ \Pi_N w_\tau + \varepsilon^2 k_3 \Pi_N w_x = \varepsilon^2 a_3 \Pi_N + \varepsilon^2 b_3 \Pi_N v + \varepsilon^2 c_3 \Pi_N w + O(\varepsilon^{N+1} e^{-\kappa\tau}), \end{cases} \quad (24)$$

where constants $C, \kappa > 0$ generally speaking depend on $N$.

## 4    Evaluation of the residual term

Theorem on the evaluation of the residual term.

*Theorem.* Let the initial conditions (2) be infinitely differentiable, and let N be any natural number. Then the solution of problem (1)-(2) can be represented as

$$\begin{pmatrix} u(x,t,\varepsilon) \\ v(x,t,\varepsilon) \\ w(x,t,\varepsilon) \end{pmatrix} = \sum_{i=0}^{N} \varepsilon^i \left( \begin{pmatrix} \bar{u}_i(x,t) \\ \bar{v}_i(x,t) \\ \bar{w}_i(x,t) \end{pmatrix} + \begin{pmatrix} \pi_i u(x,\tau) \\ \pi_i v(x,\tau) \\ \pi_i w(x,\tau) \end{pmatrix} \right) + \begin{pmatrix} R_N u(x,t,\varepsilon) \\ R_N v(x,t,\varepsilon) \\ R_N w(x,t,\varepsilon) \end{pmatrix} = \quad (25)$$
$$= U_N + \Pi_N + R_N$$

, where $R_N$ is the residual term in the expansion (25) - is the solution of the problem

$$\begin{cases} \varepsilon^2 (R_N u_t + k_1 R_N u_x) = -a R_N u + b R_N v + \varepsilon^2 c_1 R_N w + O(\varepsilon^{N+1}), \\ \varepsilon^2 (R_N v_t + k_2 R_N v_x) = a R_N u - b R_N v + \varepsilon^2 c_2 R_N w + O(\varepsilon^{N+1}), \\ R_N w_t + k_3 R_N w_x = a_3 R_N u + b_3 R_N v + c_3 R_N w + O(\varepsilon^{N+1}), \end{cases} \quad (26)$$

with initial conditions

$$\begin{cases} R_N u(x,0) = 0, \\ R_N v(x,0) = 0, \\ R_N w(x,0) = 0. \end{cases} \quad (27)$$

For the proof, it suffices to take an arbitrary natural number N, and construct N+1 terms of the expansion

$$\begin{pmatrix} u(x,t,\varepsilon) \\ v(x,t,\varepsilon) \\ w(x,t,\varepsilon) \end{pmatrix} = \sum_{i=0}^{N+1} \varepsilon^i \left( \begin{pmatrix} \bar{u}_i(x,t) \\ \bar{v}_i(x,t) \\ \bar{w}_i(x,t) \end{pmatrix} + \begin{pmatrix} \pi_i u(x,\tau) \\ \pi_i v(x,\tau) \\ \pi_i w(x,\tau) \end{pmatrix} \right) + \begin{pmatrix} R_{N+1} u(x,t,\varepsilon) \\ R_{N+1} v(x,t,\varepsilon) \\ R_{N+1} w(x,t,\varepsilon) \end{pmatrix} =$$
$$= U_{N+1} + \Pi_{N+1} + R_{N+1},$$

substitute (25) in (1) and, taking into account the algorithm for constructing the expansion terms and the estimates of (12), (24), get the estimates of (26). The initial conditions (27) for $R_N$ is obtained similarly.

*Remark.* The number of terms of the expansion N depends on the smoothness of the initial conditions (2). For infinite smoothness of the initial conditions, the number N can be arbitrary (which, of course, does not guarantee the convergence of series (4) and (13) in the usual sense).

## 5    Results

1. The main term of the AE solution of problem (1) - (2) has the form

$$\begin{pmatrix} u(x,t,\varepsilon) \\ v(x,t,\varepsilon) \\ w(x,t,\varepsilon) \end{pmatrix} = \begin{pmatrix} \bar{u}_0(x,t) \\ \bar{v}_0(x,t) \\ \bar{w}_0(x,t) \end{pmatrix} + \begin{pmatrix} \pi_0 u(x,\tau) \\ \pi_0 v(x,\tau) \\ 0 \end{pmatrix} + O(\varepsilon) =$$

$$= \begin{pmatrix} \bar{u}_0(x,t) \\ \dfrac{a}{b}\bar{u}_0(x,t) \\ \bar{w}_0(x,t) \end{pmatrix} + \begin{pmatrix} \psi(x)e^{-(a+b)\tau} \\ -\psi(x)e^{-(a+b)\tau} \\ 0 \end{pmatrix} + O(\varepsilon) \qquad (28)$$

In the domain $t > t_0$, where $t_0 > 0$ is any positive number that does not depend on $\varepsilon$, the principal term of AE (25), taking into account the exponential decrease of functions $\pi_0 u(x,\tau), \pi_0 v(x,\tau)$, takes the form

$$\begin{pmatrix} u(x,t,\varepsilon) \\ v(x,t,\varepsilon) \\ w(x,t,\varepsilon) \end{pmatrix} = \begin{pmatrix} \bar{u}_0(x,t) \\ \dfrac{a}{b}\bar{u}_0(x,t) \\ \bar{w}_0(x,t) \end{pmatrix} + O(\varepsilon) \qquad (29)$$

where $\bar{u}_0(x,t), \bar{w}_0(x,t)$ are defined as the solution of system (8) with initial conditions (20).

2. Note that the characteristics of system (8) $x - Vt = const, x - k_3 t = const$ do not coincide with the characteristics of system (1) $x - k_1 t = const, x - k_2 t = const, x - k_3 t = const$. This is not essential for smooth initial data, but it is quite essential for discontinuous initial data (2). Let one (or some of the initial data) suffer a discontinuity of the first kind at a point $x = X$. In this case, the solution of problem (1) –(2) will also have discontinuities. Well known (for example, [7],[8],[9]) that all discontinuities of the solution of system (1) can be only on characteristics, i.e. on lines $x - k_1 t = X, x - k_2 t = X, x - k_3 t = X$. However, the principal term of PAR (25), which is a solution of system (8) with discontinuous initial conditions (2), can have a discontinuity only on the characteristics of system (8), i.e. lines $x - Vt = X, x - k_3 t = X$. System (8) has two characteristics, not three, and only one of the two characteristics coincides with one of the characteristics of the original system (1). The lines $x - Vt = const$ can be called "pseudo-characteristics". This means that in asymptotically small neighborhoods of the characteristics $x - Vt = const$ of system (8), the exact solution of problem (1)-(2) remains smooth, but has the character of a narrow transition layer. In the works [4],[5],[6] and others it is shown that the solution of a singularly perturbed system of transport equations with fast variables in the critical case is described in the vicinity of the" pseudo -characteristic" - the line of discontinuity of the zero part of the regular part of the PAR, by some parabolic equation. It is of interest to see what changes in the solution behavior will result from adding slow components. In continuation of this work, we will conduct a study of AE in the case of a first-order discontinuity in the initial conditions, for example

$$\begin{cases} u(x,0,\varepsilon) = \begin{cases} u^0(x), x > 0, \\ 0, x < 0, \end{cases} \\ v(x,0,\varepsilon) = \begin{cases} v^0(x), x > 0, \\ 0, x < 0, \end{cases} \\ w(x,0,\varepsilon) = \begin{cases} w^0(x), x > 0, \\ 0, x < 0, \end{cases} \end{cases}$$

or a kind of "narrow cap", "splash"

$$\begin{cases} u(x,0,\varepsilon) = u^0(x/\varepsilon), \\ v(x,0,\varepsilon) = v^0(x/\varepsilon), \\ w(x,0,\varepsilon) = w^0(x/\varepsilon), \end{cases}$$

where $|u^0(z)|, |v^0(z)|, |w^0(z)| < Ce^{-\kappa z^2}, C, \kappa > 0.$ In this case AE of the solution of the problem (1)-(2) will have a more complex form.